\documentclass[11pt,a4paper]{article}

\usepackage[a4paper,margin=1.9cm]{geometry}
\usepackage{xcolor}
\usepackage{graphicx}
\usepackage{epstopdf}
\usepackage{tabularx}
\usepackage{booktabs}
\usepackage{array}
\usepackage{enumitem}
\usepackage{titlesec}
\usepackage{fancyhdr}
\usepackage{hyperref}
\usepackage{needspace}
\usepackage[most]{tcolorbox}
\usepackage{subcaption}
\usepackage{amssymb}

\definecolor{deepblue}{HTML}{2D5A27}
\definecolor{teal}{HTML}{4A7C2F}
\definecolor{leaf}{HTML}{5B8C5A}
\definecolor{gold}{HTML}{C18B2A}
\definecolor{softblue}{HTML}{EEF5E8}
\definecolor{softgreen}{HTML}{E6F2E6}
\definecolor{softgold}{HTML}{FFF6E3}
\definecolor{linegray}{HTML}{D8DEE4}
\definecolor{textgray}{HTML}{3A3A3A}
\definecolor{brightblue}{HTML}{1A73E8}

\hypersetup{
  colorlinks=true,
  linkcolor=deepblue,
  urlcolor=brightblue,
  citecolor=brightblue
}

\renewcommand{\arraystretch}{1.22}

\titleformat{\section}
  {\Large\bfseries\color{deepblue}}
  {}
  {0pt}
  {}
\titleformat{\subsection}
  {\large\bfseries\color{deepblue}}
  {}
  {0pt}
  {}

\renewcommand{\headrulewidth}{0.3pt}
\renewcommand{\headrule}{\hbox to\headwidth{\color{linegray}\leaders\hrule height \headrulewidth\hfill}}

\newcolumntype{Y}{>{\raggedright\arraybackslash}X}
\newcolumntype{L}[1]{>{\raggedright\arraybackslash}p{#1}}
\newcolumntype{R}{>{\raggedleft\arraybackslash}p{0.12\textwidth}}

\newtcolorbox{briefbox}[2][]{
  enhanced,
  breakable,
  colback=#2,
  colframe=#2!70!black,
  boxrule=0.6pt,
  arc=1.5mm,
  left=3mm,
  right=3mm,
  top=2mm,
  bottom=2mm,
  #1
}

\newtcolorbox{metricbox}[2][]{
  enhanced,
  colback=#2,
  colframe=#2!70!black,
  boxrule=0.6pt,
  arc=1.5mm,
  left=3mm,
  right=3mm,
  top=2mm,
  bottom=2mm,
  #1
}

\newcommand{\metric}[2]{%
  \Needspace{8\baselineskip}%
  \begin{metricbox}{softblue}
  {\large\bfseries\color{deepblue}#1}\par
  #2
  \end{metricbox}
}

\newcommand{\sumlablogo}{%
  \IfFileExists{sum_logo.png}{%
    \includegraphics[width=3.4cm]{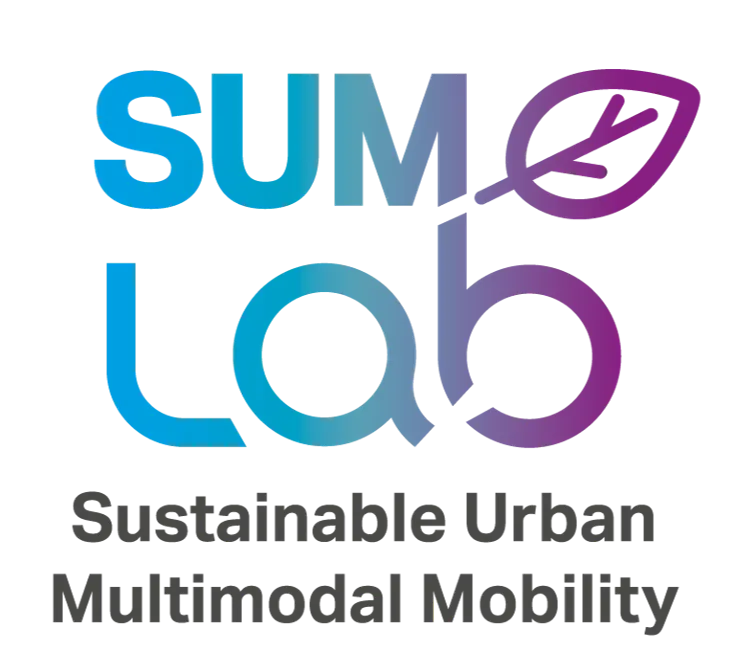}%
  }{%
    \begin{minipage}{3.4cm}
      \raggedleft
      {\Large\bfseries\color{teal}SUM\par}
      {\Large\bfseries\color{deepblue}Lab\par}
      {\scriptsize\color{textgray}Sustainable Urban\par Multimodal Mobility\par}
    \end{minipage}%
  }%
}

\newcommand{\keymessage}[3]{%
  \begin{tcolorbox}[
    enhanced,
    colback=white,
    colframe=linegray,
    boxrule=0.5pt,
    arc=1.5mm,
    left=1.3cm,right=3mm,top=2.5mm,bottom=2.5mm,
    before skip=0.45em,after skip=0.1em,
    overlay={%
      \fill[deepblue,rounded corners=1.2mm]
        ([xshift=0.5pt,yshift=-0.5pt]frame.north west) rectangle
        ([xshift=1.2cm,yshift=0.5pt]frame.south west);
      \node[text=white,font=\LARGE\bfseries,inner sep=0pt] at
        ([xshift=0.6cm]frame.west) {#1};%
    }
  ]
  \textbf{\color{deepblue}#2}\par\smallskip
  \small #3
  \end{tcolorbox}%
}

\newcommand{\audheader}[1]{%
  \begin{tcolorbox}[
    enhanced,
    colback=deepblue,
    colframe=deepblue,
    boxrule=0pt,
    arc=1.5mm,
    left=3mm,right=3mm,top=2.5mm,bottom=2.5mm,
    before skip=0.9em,after skip=0.5em,
  ]
  {\small\bfseries\color{white}\MakeUppercase{#1}}%
  \end{tcolorbox}%
}

\newcommand{\recitem}[2]{%
  \begin{tcolorbox}[
    enhanced,
    colback=white,
    colframe=teal,
    leftrule=3pt,toprule=0pt,bottomrule=0pt,rightrule=0pt,
    arc=0pt,outer arc=0pt,
    left=3mm,right=3mm,top=2mm,bottom=2mm,
    before skip=0.5em,after skip=0.1em,
  ]
  {\small$\blacktriangleright$~\textbf{#1}}\par\smallskip
  {\small #2}
  \end{tcolorbox}%
}

\newcommand{\calloutbox}[2]{%
  \begin{tcolorbox}[
    enhanced,
    colback=softblue,
    colframe=teal,
    leftrule=3pt,toprule=0pt,bottomrule=0pt,rightrule=0pt,
    arc=0pt,outer arc=0pt,
    left=3mm,right=3mm,top=2mm,bottom=2mm,
    before skip=0.5em,after skip=0.2em,
  ]
  {\scriptsize\bfseries\color{teal}\MakeUppercase{#1}}\par\smallskip
  {\small #2}
  \end{tcolorbox}%
}

\begin{document}

\hypersetup{pageanchor=false}
\begin{titlepage}
\pagecolor{softblue}
\vspace*{0.30cm}
\begin{minipage}[t]{0.66\textwidth}
{\Large\bfseries\color{teal} STRATEGIC POLICY BRIEF - 2026\par}
\end{minipage}
\hfill
\begin{minipage}[t]{0.28\textwidth}
\raggedleft
\vspace{-0.20cm}
\sumlablogo
\end{minipage}
\vspace{0.45cm}

{\Huge\bfseries\color{deepblue} Modular Autonomous Transit\par}
\vspace{0.22cm}
{\LARGE\bfseries\color{deepblue} From Vehicle Modularity\par}
\vspace{0.08cm}
{\LARGE\bfseries\color{deepblue} to Deployment-Ready Transit Operations\par}
\vspace{0.18cm}
{\color{teal}\rule{1.7cm}{1.0pt}}\hspace{0.12cm}{\color{leaf}\rule{1.15cm}{1.0pt}}\hspace{0.12cm}{\color{gold}\rule{0.85cm}{1.0pt}}\par
\vspace{0.26cm}
{\large An uncertainty-aware five-layer deployment roadmap for modular autonomous vehicles\par}
\vspace{0.26cm}
\begin{center}
{\small\bfseries\color{teal} PREPARED BY\par}
\vspace{0.06cm}
{\LARGE\bfseries\color{deepblue} Dongyang Xia\par}
\vspace{0.04cm}
{\LARGE\bfseries\color{deepblue} Shadi Sharif Azadeh\par}
\vspace{0.04cm}
{\normalsize Sustainable Urban Multi-modal Mobility (SUM) Lab, Delft University of Technology\par}
\end{center}
\vspace{0.16cm}
\begin{briefbox}[colframe=linegray, boxrule=0.5pt, left=4mm, right=4mm, top=2.5mm, bottom=2.5mm]{white}
\small
\textbf{Audience:} Tech companies such as NExT leadership; public transport operators; autonomous transit developers; charging partners; public authorities, and industrial partners considering pilots.

\textbf{Content:} The strongest aspect of Modular Autonomous Vehicles (MAVs) is not only that vehicles can physically couple. A MAV is composed of Modular Autonomous Units (MAUs), and modularity creates operating advantages across line-level operations, intermodal services, network-level operations, stochastic real-time control, and electrified vehicle operations and charging infrastructure.

\smallskip
{\bfseries\color{deepblue} Five uncertainty-aware operating layers}

\begingroup
\scriptsize
\renewcommand{\arraystretch}{1.08}
\begin{tabularx}{\linewidth}{@{}L{0.32\linewidth}Y@{}}
\toprule
\textbf{Layer} & \textbf{Research signal} \\
\midrule
Line-level operations & Results based on Beijing Bus Line 468 indicate that station-wise coupling and decoupling can reduce used MAUs and operating cost under uncertain demand. \\
\addlinespace
Intermodal services & Results indicate that one fleet of MAUs can serve both fixed-line and demand-responsive services under stochastic demand. \\
\addlinespace
Network-level operations & Experiments based on a Beijing bus subnetwork indicate the value of cross-line circulation of MAUs, in-vehicle transfers of passengers, and integrated optimization of timetabling and vehicle scheduling under uncertainty. \\
\addlinespace
Stochastic real-time control & Results based on a Beijing bus subnetwork indicate that timetabling and vehicle scheduling of MAUs can be rescheduled as real-time demand information knows within seconds. \\
\addlinespace
Electrified vehicle operations and charging infrastructure & Results using data from NExT MAUs indicate that entire-day electric operation needs vehicle scheduling, charging, and charger deployment. \\
\bottomrule
\end{tabularx}
\endgroup

\smallskip
\textbf{Research basis and disclaimer:} the brief summarizes scientific findings of four studies conducted by Dr. Xia and Dr. Sharif Azadeh on modular transit operations, including a charging infrastructure study using technical data associated with NExT MAUs.
\normalsize
\end{briefbox}
\vfill
\noindent\makebox[\textwidth][c]{\small\color{textgray}\today}\par
\end{titlepage}
\hypersetup{pageanchor=true}
\pagecolor{white}

\section{Key Messages}

\keymessage{1}{Use a five-layer deployment framework.}{The research supports a clear progression of line-level operations, intermodal services, network-level operations, stochastic real-time control, and electrified vehicle operations and charging infrastructure. The first three layers optimize modular operations under demand uncertainty; the fourth layer moves uncertainty into real-time control; the fifth layer embeds vehicle scheduling in charging station location, charger deployment, and all-day electric operation.}

\keymessage{2}{Line-level operations: Station-wise coupling and decoupling reduce wasted capacity.}{In the Beijing Bus Line 468 case study, allowing MAVs to couple and decouple along the line reduced the number of used MAUs by about 46\% and operating costs by about 19\%, compared with a strategy that only changes capacity over time. Passengers' waiting costs increased by only about 0.8\%.}

\keymessage{3}{Intermodal services: The same fleet can support fixed-line and demand-responsive services.}{Integrated fixed-line and station-to-station demand-responsive transit setting reduced the combined passenger and operating costs by about 6\% compared with a practical uniform-headway, fixed-capacity bus operation.}

\keymessage{4}{Network-level operations: Cross-line circulation is where modularity becomes strategically important.}{In a Beijing bus subnetwork with two bidirectional lines, 89 stations, up to 50 trips, and a four-hour operating horizon, fully flexible modular operation used fewer MAUs across the tested instances. In one large instance, fixed-capacity operation required 96 MAUs, while the fully flexible strategy required 69.}

\keymessage{5}{Stochastic real-time control: Modularity creates operating resilience.}{The network study shows that a base timetable and vehicle schedule can be fine-tuned as demand information arrives. The learning-based real-time framework produced high-quality decisions within one minute per decision stage, and a 10\% temporary overload allowance absorbed a 5\% demand increase without additional capacity resources.}

\keymessage{6}{Electrified vehicle operations and charging infrastructure: Electric MAUs cannot be planned with depot charging alone.}{Charging infrastructure and operations research based on data from NExT MAVs shows that vehicle scheduling, charging-station location, charger deployment, depot charging, and en-route fast charging need to be considered together. In an all-day case from 07:00 to 22:00, the optimized plan selected two en-route fast-charging stations, each with two fast chargers, and used a fleet of 101 MAUs.}

\keymessage{7}{Manufacturers need to deliver vehicles, interfaces, and operational data.}{Dynamic modularity becomes valuable only when demand forecasts, real-time loads, vehicle locations, formation status, battery state, and charger occupancy are available to the dispatch system.}




\newpage

\section{Recommendations}

\phantomsection\label{sec:nextrec}%
\audheader{For NExT and Modular Autonomous Vehicle Manufacturers}

\recitem{Sell modularity as a data-driven deployment capability.}{The commercial argument should move from ``the vehicles can couple'' to ``under this demand, network, station, and charging context, the system reduces idle capacity, fleet units, waiting time, and operating cost.'' System economics, not vehicle novelty, is the strongest opening for operator and city discussions.}

\phantomsection\label{rec:M7}%
\recitem{Offer a pilot and visualization platform.}{A strong offer includes candidate corridors, station requirements, docking and passenger-circulation rules, charging infrastructure analysis, dispatch APIs, before-after KPIs, and a dashboard that makes modularity legible to executives, operators, regulators, and infrastructure partners. Recommended metrics include the number of used MAUs, cost per unit of capacity, average waiting time, load factor, unserved passengers, charging-queue time, and failure-recovery time. A research partner such as SUM Lab can support this by building demand scenarios, optimization models, solution algorithms, and a visualization platform.}

\recitem{Deliver data interfaces.}{Formation state, passenger load, state of charge, charger occupancy, docking availability, task state, and failure mode are prerequisites for scheduling value. A vehicle delivered with these interfaces can be integrated into dispatch algorithms instead of being evaluated only as a standalone vehicle.}


\recitem{Package the vehicle with entire-day electric feasibility.}{The charging study using NExT technical specifications supports a credible claim that NExT MAUs can be planned for entire-day operations when charging station location, charger deployment, partial charging, en-route fast charging, and unit schedules are optimized together.}

\phantomsection\label{sec:oprec}%
\audheader{For Public Transport Operators}

\recitem{Do not evaluate MAVs only as replacements for conventional buses.}{A pilot should explicitly state whether it allows station-wise coupling, cross-line scheduling, demand-responsive transit tasks, and en-route fast charging.}

\recitem{Evaluate passengers' costs and operating costs together.}{The research shows that the trade-off is manageable; optimizing only one side can understate the value of MAV fleets.}

\recitem{Move from a single line toward a small network pilot.}{A single line can demonstrate capacity matching. A small network can demonstrate cross-line circulation, in-vehicle transfer, demand-responsive tasks, and fleet-size savings.}

\recitem{Use auditable service metrics from the start.}{Pilot evaluation should include passengers' waiting cost, average waiting time, unserved passengers, load factor, fleet size, charging-queue time, service reliability, and recovery time after vehicle or charger disruptions.}

\newpage 
\phantomsection\label{sec:citrec}%
\audheader{For Cities and Regulators}

\recitem{Reserve space for coupling, decoupling, and en-route fast charging.}{Modular transit needs more than road-testing permission; it needs workable station and charging layouts.}

\recitem{Establish operational data-sharing rules.}{Cities need demand, load factors, vehicle kilometers, charging needs, and service reliability to evaluate whether modular transit improves public transport.}

\recitem{Use pilot zones to validate safety and service standards.}{Key issues include docking safety, passenger passage between units, emergency evacuation, degraded formation after failures, manual takeover procedures, and operating-data transparency.}

\audheader{For Charging Infrastructure Partners and Energy Companies}

\recitem{Treat charging as an operational feasibility problem.}{Depot charging alone is insufficient for all-day electric operation in the studied setting. Charging partners should position en-route fast charging as a functional requirement for continuous service, not as an optional add-on.}

\recitem{Co-plan chargers with vehicle schedules.}{Fast-charging station location, charger quantity, partial-charging rules, grid capacity, charger reliability, and queue management should be planned jointly with timetables, formation decisions, and MAU schedules.}

\recitem{Provide charger data interfaces.}{Charger availability, occupancy, power limits, fault status, queue length, and energy price signals should be visible to the dispatch system so that charging decisions can be optimized during service.}

\newpage

\section{Five-Layer Research Findings}

\begin{briefbox}{softgreen}
The studies are best read as a deployment ladder for transit operators, technology developers, and public authorities evaluating MAV pilots. The first three layers address operations under demand uncertainty at increasing scale, from single-line to intermodal and network-wide. The fourth layer extends this to real-time scheduling as demand evolves during service. The fifth layer integrates vehicle operations with charging station location, charger deployment, and all-day electric operation.

\smallskip
Modularity's value depends on how the vehicles are operated. If units are coupled into fixed formations and dispatched like conventional buses, most of the fleet and cost advantage disappears. All five layers rely on three enabling capabilities working in coordination. \textbf{Vehicle capability} means units couple and decouple reliably with low delay. \textbf{Dispatch capability} means timetables, schedules, and charging are optimized jointly. \textbf{Data capability} means the system knows where demand is, where units are, and which units need charging.
\end{briefbox}

\begin{tabularx}{\textwidth}{L{0.28\textwidth}YL{0.22\textwidth}}
\toprule
\textbf{Layer} & \textbf{Quantified signal} & \textbf{Study} \\
\midrule
Line-level operations & \textbf{-46\%} used MAUs and \textbf{-19\%} operating cost; waiting cost: \textbf{+0.8\%} (baseline: buses that fix their size before each trip and do not resize at stops; Beijing Bus Line 468). & \href{https://doi.org/10.1016/j.trc.2023.104314}{Xia et al. (2023)}, \textit{TR-C} \\
\addlinespace
Intermodal services & \textbf{-6\%} combined passenger and operating cost (baseline: standard fixed-route, fixed-size buses; Beijing Bus Line 468). & \href{https://doi.org/10.1016/j.trc.2024.104610}{Xia et al. (2024)}, \textit{TR-C} \\
\addlinespace
Network-level operations & Fleet size reduced by \textbf{28\%--51\%} across operational scenarios (baseline: buses running at fixed vehicle sizes; Beijing bus subnetwork). & \href{https://doi.org/10.1287/trsc.2025.0116}{Xia et al. (2026)}, \textit{Transportation Science} \\
\addlinespace
Stochastic real-time control & Joint real-time timetable and vehicle-schedule reoptimization cuts total cost by \textbf{3.15\%--8.14\%} vs.\ a fixed plan; vehicle-only scheduling saves \textbf{0.81\%--3.47\%} (baseline: no reoptimization; Beijing bus subnetwork). & \href{https://doi.org/10.1287/trsc.2025.0116}{Xia et al. (2026)}, \textit{Transportation Science}  \\
\addlinespace
Electrified vehicle operations and charging infrastructure & Locking vehicle size at 2 units raises passengers' cost by \textbf{+16.22\%}; locking at 4 units by \textbf{+28.22\%} (baseline: freely resizable electric MAUs; NExT vehicle data). & \href{https://arxiv.org/abs/2504.04408}{Xia et al. (2025)} \textit{TR-B} (under review) \\
\bottomrule
\end{tabularx}

\section{Management Insights}

\subsection{1. Line Level Operations: Station-wise Coupling and Decoupling Turn Mismatch into a Manageable Problem}

The first study (\href{https://doi.org/10.1016/j.trc.2023.104314}{Xia et al., 2023}) uses operational data from Beijing Bus Line 468 to optimize timetables and dynamic capacity allocation for MAVs on a single fixed-route line. The model reflects two practical realities. Passenger demand is uneven across time and stations, and future demand is uncertain.

\metric{Insights}{Passenger demand is uneven across both time and space. Flexible capacity that adjusts only over time leaves the station-wise imbalance unresolved. Adding station-wise coupling and decoupling on top of time-varying capacity further reduced the number of utilized MAUs by about 46\% and operating costs by about 19\%, compared with time-varying capacity alone, with only a 0.8\% increase in passengers' waiting costs.}

The most valuable use of MAVs is to place capacity exactly where it is needed. Some stations and sections need more MAUs during the peak. Low demand sections do not need a large formation to continue running with empty capacity. Fixed capacity operation can create both waste and crowding, with too much capacity in low-demand trips and too little in high-demand trips.

\calloutbox{Manufacturer action}{station-wise coupling and decoupling is not yet a tested capability in revenue service. Manufacturers need to demonstrate docking time under realistic station conditions, safety envelopes that fit within existing station layouts, passenger movement between coupled units, and service continuity after individual unit failures. Without this data, the optimization value shown in the research cannot be translated into a procurement decision.}

\calloutbox{Data required}{station-wise boarding and alighting, time-dependent running times, arrival and departure times, load factors, and real-time unit availability.}

\calloutbox{Cost of action}{The main investment is not necessarily more vehicles. It is the integration of vehicle-control systems with dispatch systems, plus operational validation at selected stations.}

\Needspace{15\baselineskip}
\subsection{2. Intermodal Services: One Modular Fleet Serves Both Fixed-Line and Demand-Responsive Transit}

The second study (\href{https://doi.org/10.1016/j.trc.2024.104610}{Xia et al., 2024}) integrates fixed-line scheduled service with station-to-station demand-responsive transit. Modular units can be assigned to fixed-route trips, dispatched for demand-responsive requests, and then returned to fixed-line operations.

\metric{Insights}{The same pool of MAVs simultaneously served scheduled fixed-line trips and on-demand station-to-station requests during the same operating period, without a separate demand-responsive fleet. This dual use, optimized jointly, reduced expected combined passenger and operating costs by about 6\% compared with a conventional uniform-headway, fixed-capacity fixed-line service.}

This matters commercially because it expands the market position of MAVs. They do not have to be sold only as replacements for conventional buses. They can be sold as a platform for multimodal public transit, combining fixed-route service for stable corridors and demand-responsive service for low-density or scattered demand.

\calloutbox{Manufacturer action}{Support service-mode switching at the vehicle and back-office level, covering fixed-line service, short-turning service, demand-responsive service, deadheading, standby, and charging. The dispatch system should know whether each unit is available, dockable, passenger-ready, and compatible with its next task.}

\calloutbox{Data required}{Fixed-line origin-destination demand, demand-responsive requests, waiting times, unit task status, feasible pickup and drop-off stations, and passenger tolerance for transfers or detours.}

\calloutbox{Cost of action}{Operators need to reorganize fleet control, but they do not necessarily need a separate demand-responsive fleet. The commercial message is that one asset base can serve two types of demand.}

\subsection{3. Network-Level Operations: Cross-Line Circulation and In-Vehicle Transfers Are Scale-Up Features}\label{sec:layer3}

The third study (\href{https://doi.org/10.1287/trsc.2025.0116}{Xia et al., 2026}) moves from one line to a network. Modular units can circulate across lines and change formations at transfer stations. Passengers can also transfer inside the vehicle by walking from one coupled unit to another, without alighting and reboarding.

\metric{Insights}{For the three largest tested network instances, sequential planning (designing the timetable first under assumed unlimited resources, then scheduling vehicles) found no feasible solution. The resulting timetable could not be executed with the available fleet. Integrated planning solved all instances. Cross-line unit circulation then reduced fleet requirements by 28--51\% compared with fixed-composition operation across the same network.}

The numerical pattern is consistent across different instance sizes. The network experiments use a Beijing bus subnetwork with two bidirectional lines, 89 stations, up to 50 trips, and a four-hour operating horizon. The fully flexible strategy usually reduces operating costs substantially, while keeping passengers' costs similar to the fixed-capacity strategy.

{\small
\begin{tabularx}{\textwidth}{p{0.18\textwidth}p{0.13\textwidth}p{0.13\textwidth}p{0.13\textwidth}Y}
\toprule
\textbf{Instance} & \textbf{Fixed units} & \textbf{Flexible units} & \textbf{Fleet reduction} & \textbf{Service impact} \\
\midrule
Smaller instances & 72 & 35 & \textbf{-51.4\%} & Passengers' cost slightly lower; operating cost reduced by about 46\%. \\
\addlinespace
Medium instances & 96 & 51 & \textbf{-46.9\%} & Passengers' cost slightly lower; operating cost reduced by about 44\%. \\
\addlinespace
Larger instances & 96 & 69 & \textbf{-28.1\%} & Passengers' cost nearly unchanged; operating cost reduced by about 36\%. \\
\bottomrule
\end{tabularx}
}

This finding shifts MAVs from a line-efficiency tool to a network-design tool. Traditional sequential planning first designs the timetable and then schedules vehicles. In larger modular networks, that approach can produce timetables that cannot be operated with the available fleet. Integrated planning avoids this by designing the timetable and vehicle schedule together.

{\small
\begin{tabularx}{\textwidth}{p{0.20\textwidth}p{0.16\textwidth}p{0.16\textwidth}Y}
\toprule
\textbf{Instance} & \textbf{Sequential} & \textbf{Integrated} & \textbf{Operational interpretation} \\
\midrule
Smallest instances & 21 units & 21 units & Similar performance at the smallest scale. \\
\addlinespace
Smaller instances & 36 units & 35 units & Integration begins to save fleet resources. \\
\addlinespace
Medium instances & 59 units & 51 units & Integration reduces required units by 8. \\
\addlinespace
Larger instances & No feasible solution & 58, 61, and 69 units & Sequential planning fails at network scale; integrated planning solves all cases. \\
\bottomrule
\end{tabularx}
}

\begin{figure}[h]
     \centering
     \begin{subfigure}{0.45\textwidth}
         \centering
         \includegraphics[width=\textwidth]{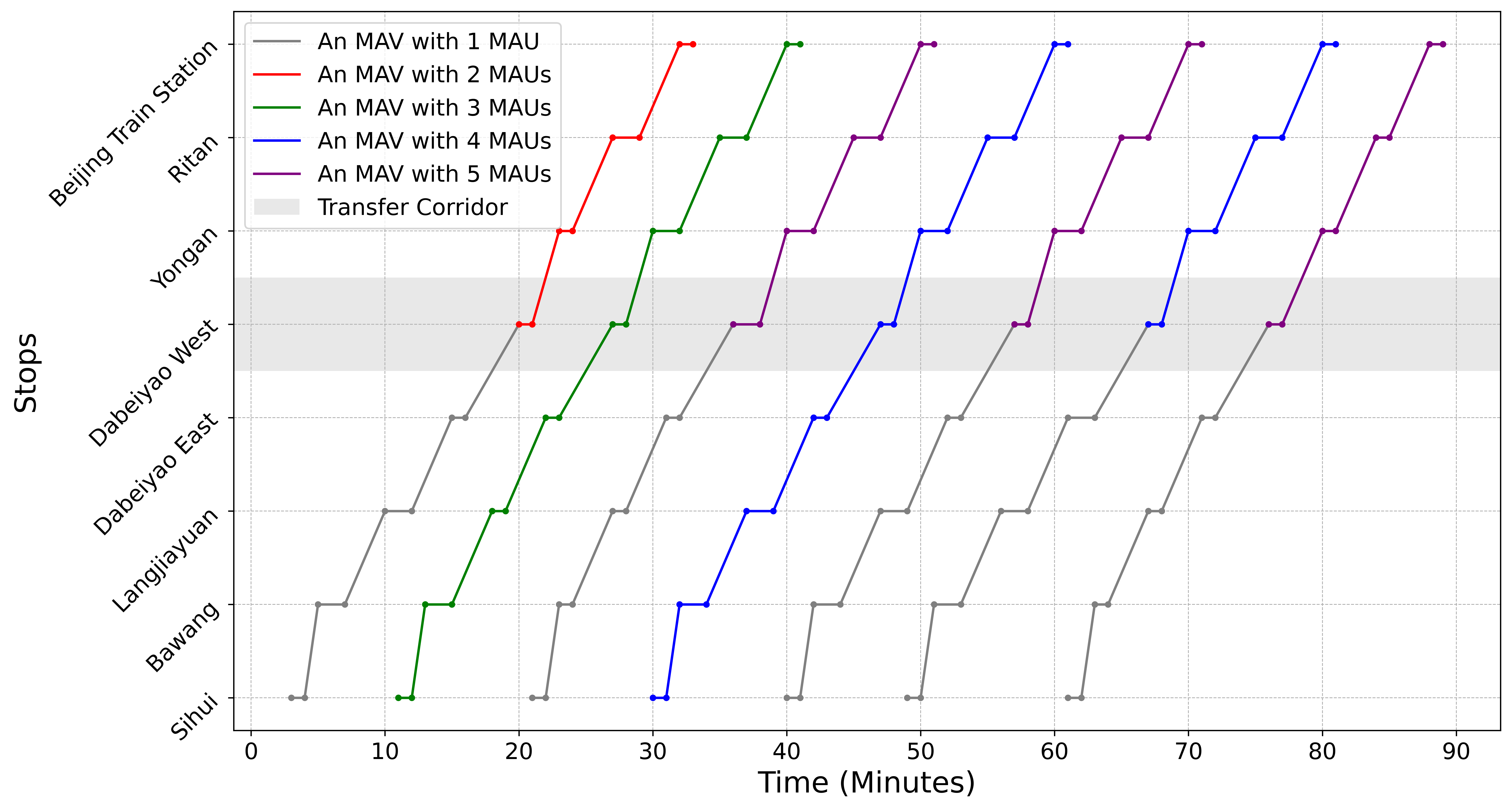}
         \caption{The timetable and formations of MAVs}
     \end{subfigure}
     \begin{subfigure}{0.49\textwidth}
         \centering
         \includegraphics[width=\textwidth]{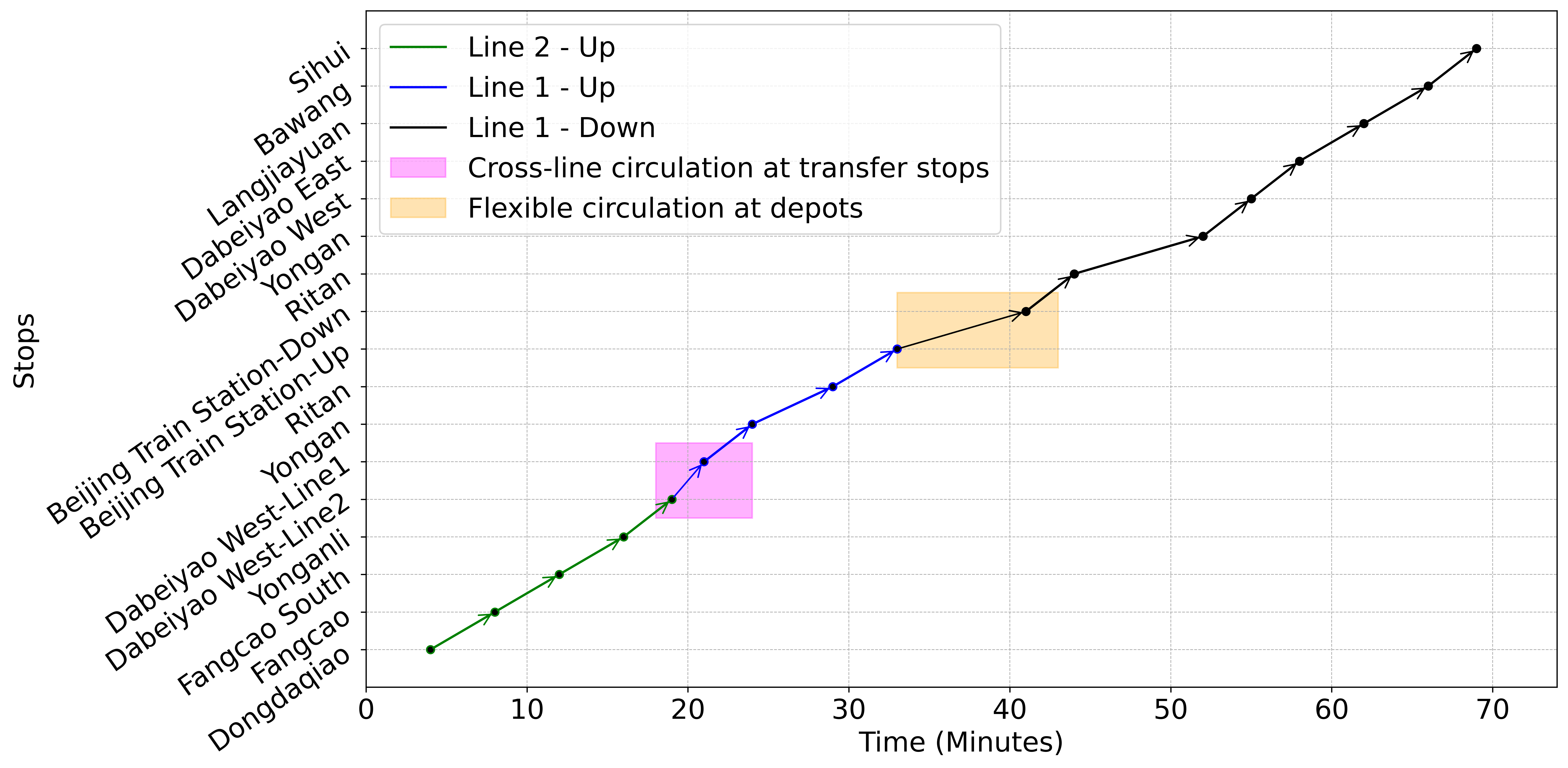}
         \caption{Movements of an MAU}
     \end{subfigure}
        \caption{Optimized timetable and unit formations for a two-line network, showing how MAVs couple and decouple across stations and circulate between lines. The shaded station is the transfer hub where in-vehicle passenger transfers occur without alighting. Source: Xia et al.\ (2026), \textit{Transportation Science}.}
        \label{fig:resultTTVS}
\end{figure}

The benefits for passengers are also important. In-vehicle transfers reduce the inconvenience of getting off, walking, waiting, and boarding again. For manufacturers, this turns open gangways, passenger circulation, docking precision at transfer stations, and vehicle-to-vehicle communication into product features with direct service value.

\calloutbox{Manufacturer action}{Treat cross-line circulation and in-vehicle transfer as advanced deployment scenarios, not as afterthoughts. A strong pilot should include at least one transfer hub, two lines, and multiple demand periods.}

\calloutbox{Data required}{Multi-line timetables, transfer origin-destination demand, station-space constraints, interline repositioning times, real-time unit positions, and operating rules for cross-line service.}

\calloutbox{Cost of action}{Operators and manufacturers need to redesign dispatch workflows together. A practical starting point is a limited set of transfer stations before expanding to network-level control.}

\subsection{4. Stochastic Real-Time Control: Modular Transit Needs Continuous Adjustment}\label{sec:layer4}

The same study (\href{https://doi.org/10.1287/trsc.2025.0116}{Xia et al., 2026}) also develops a real-time decision-making framework. Historical demand scenarios are used to produce a base timetable and vehicle schedule. During operation, the system fine-tunes departures and reoptimizes unit schedules as new demand information arrives.

\metric{Insights}{A 10\% temporary overload allowance (slightly relaxing load limits on individual sections) absorbed a 5\% demand surge without additional units or unserved passengers. When demand fell, the base plan remained feasible throughout. The real-time reoptimization framework, using machine-learning scenario selection, produced updated schedules within one minute per decision stage.}

This layer is not about optimizing the vehicle technology itself. It is about optimizing how MAVs are operated in a stochastic environment. The robustness tests translate uncertainty into operational tolerances, specifically how much demand variation the system can absorb before extra capacity or real-time intervention becomes necessary.

{\small
\begin{tabularx}{\textwidth}{p{0.23\textwidth}p{0.23\textwidth}p{0.20\textwidth}Y}
\toprule
\textbf{Demand scenario} & \textbf{Expected overloaded or unserved passengers} & \textbf{Overloaded traversals} & \textbf{Interpretation} \\
\midrule
Demand decreases of 5--20\% & 0 & 0\% & Base plan remains feasible across all tested demand decreases. \\
5\% increase & 147 & 0.72\% & Capacity stress appears but remains localized. \\
5\% increase with 10\% overload allowance & 0 & 0\% & Small operating flexibility can absorb the demand shock. \\
\bottomrule
\end{tabularx}
}

The product implication is clear. Modular vehicles should be dispatchable intelligent assets, not passive vehicles assigned to fixed trips. Vehicle-side systems need to continuously report state, and dispatch-side systems need to reassign formations, lines, and tasks quickly.

\calloutbox{Manufacturer action}{Provide real-time status interfaces and an operational dashboard. Useful fields include formation status, load factor, battery state, docking availability, failure mode, and compatibility with the next assigned task.}

\calloutbox{Data required}{Real-time demand, vehicle position, passenger load, station crowding, disruptions, task-completion times, and battery state.}

\calloutbox{Cost of action}{The main costs are software integration and data governance. The more standardized the data interface and operating rules, the easier it becomes to turn real-time control into operating value.}

\subsection{5. Electrified Vehicle Operations and Charging Infrastructure}\label{sec:layer5}

The fourth study (\href{https://arxiv.org/abs/2504.04408}{Xia et al., 2025}) focuses on electrified vehicle operations and charging infrastructure for NExT MAUs. It jointly decides where to locate fast-charging stations, how many chargers to deploy, how many units to purchase, which units operate regular or short-turning services, and when and where each unit charges.

\metric{Insights}{Even with a 90\,kWh battery and 600\,kW fast-charging capability, all-day operation (07:00--22:00) requires en-route fast charging, as depot charging alone does not provide sufficient coverage for continuous duty cycles. Under the same fleet investment, locking units to a fixed two-unit composition increased passengers' costs by 16.22\%; locking to four units increased them by 28.22\%.}

This is the study most directly linked to data from NExT MAVs in the brief. The experiments use real-life bus-line data and technical specifications of NExT MAVs, including NExT-provided inputs on battery and energy performance. They therefore support a concrete message for NExT. Electrification should be marketed as an integrated vehicle, charging infrastructure, charger deployment, and scheduling system, not as a standalone range claim.

Flexible composition also outperformed fixed composition under comparable fleet investment. Fixed two-unit electric MAVs increased passengers' costs by 16.22\%, while fixed four-unit electric MAVs increased passengers' costs by 28.22\%. Average waiting times also worsened under fixed compositions.

{\small
\begin{tabularx}{\textwidth}{p{0.25\textwidth}p{0.20\textwidth}p{0.14\textwidth}Y}
\toprule
\textbf{Comparison setting} & \textbf{Composition strategy} & \textbf{Fleet investment} & \textbf{Service impact} \\
\midrule
Flexible MAVs & Dynamic composition & 94 units & Baseline passengers' cost; average waiting time: 5.72 minutes. \\
\addlinespace
Comparable fixed fleet & Fixed two-unit electric MAVs & 94 units & Passengers' cost: +16.22\%; average waiting time: +17.46\%. \\
\addlinespace
Comparable fixed fleet & Fixed four-unit electric MAVs & 96 units & Passengers' cost: +28.22\%; average waiting time: +29.14\%. \\
\bottomrule
\end{tabularx}
}

\begin{figure}[h]
    \centering
    \includegraphics[width=0.9\linewidth]{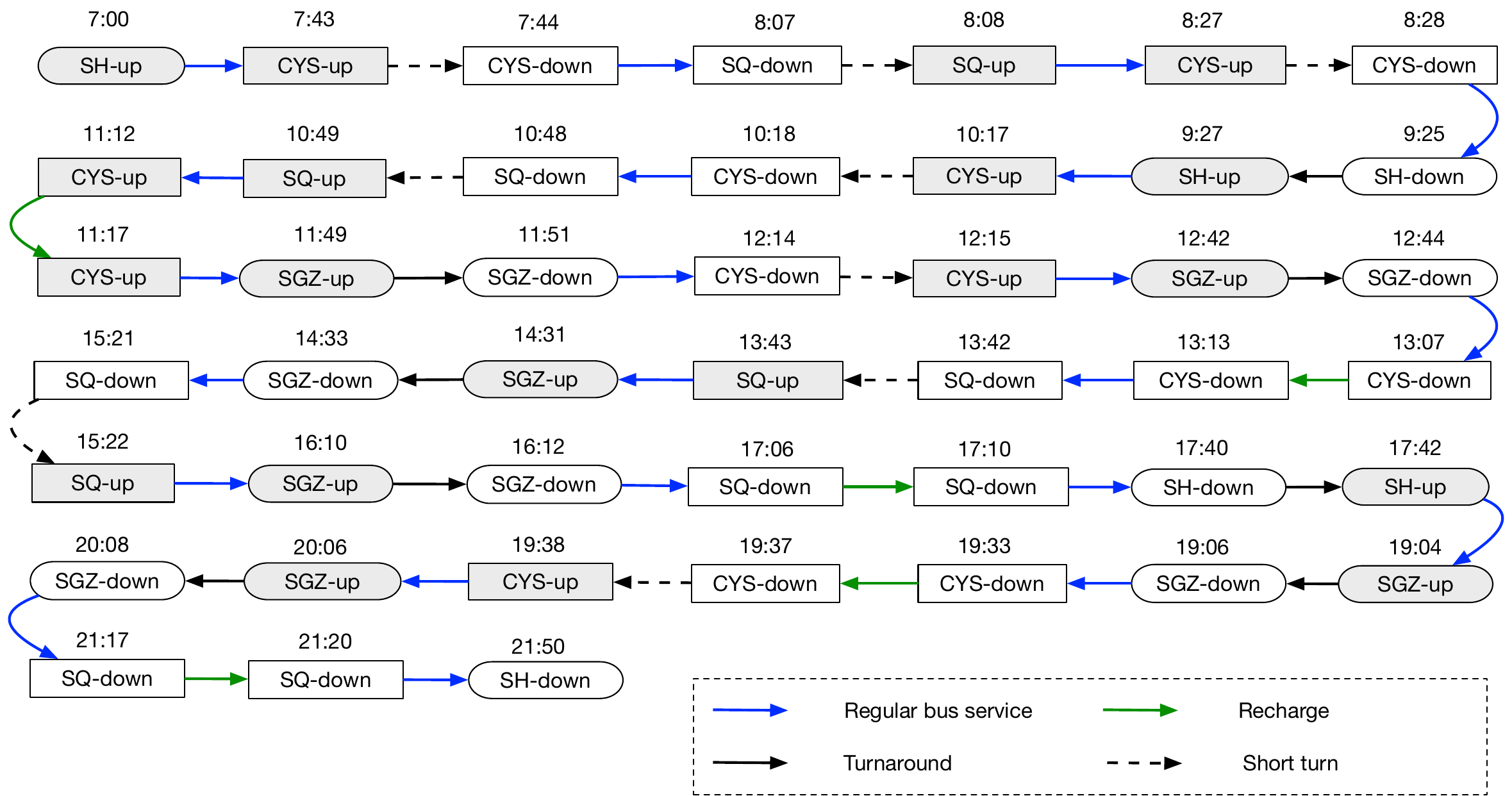}
    \caption{Optimized circulation schedule of an E-MAU between depots and charging stations from 7:00 to 22:00 with forward-only movements at non-charging stations.}
    \label{fig:optimizedCirculationAfter}
\end{figure}

Figure \ref{fig:optimizedCirculationAfter} shows the optimized circulation schedule of an E-MAU from 7:00 to 22:00, reporting the schedules between depots and charging stations. This schedule indicates that the unit undergoes multiple recharging, decoupling, and coupling operations at charging stations. Additionally, it performs a combination of regular bus and short-turning services on an intermittent basis. For instance, the unit departs from charging station SH in the up-operational direction to perform a regular bus service, arrives at station CYS, performs a short-turning service here, and returns to station SQ in the down-operational direction.

\calloutbox{Manufacturer action}{Sell the vehicle, battery, fast-charging strategy, partial-charging logic, and operations planning as one package. The key commercial question is not maximum range in isolation; it is whether the fleet can complete an all-day duty plan under uncertain demand and constrained chargers.}

\calloutbox{Data required}{Energy-consumption curves, state of charge, charging power, charger occupancy, station service constraints, vehicle duty chains, and service-reliability indicators.}

\calloutbox{Cost of action}{Charging infrastructure requires coordination and investment, but well-placed en-route fast charging can reduce dependence on oversized depot charging and excessive spare vehicles.}

\section{Research References}

\begin{enumerate}[leftmargin=1.0cm,label={[\arabic*]}]
  \item Dongyang Xia, Jihui Ma, Shadi Sharif Azadeh, and Wenyi Zhang (2023). Data-driven distributionally robust timetabling and dynamic-capacity allocation for automated bus systems with modular vehicles. \textit{Transportation Research Part C}, 155, 104314. \url{https://doi.org/10.1016/j.trc.2023.104314}

  \item Dongyang Xia, Jihui Ma, and Shadi Sharif Azadeh (2024). Integrated timetabling and vehicle scheduling of an intermodal urban transit network: A distributionally robust optimization approach. \textit{Transportation Research Part C}, 162, 104610. \url{https://doi.org/10.1016/j.trc.2024.104610}

  \item Dongyang Xia, Jihui Ma, and Shadi Sharif Azadeh (2026). Integrated timetabling and scheduling of modular autonomous vehicles under uncertainty. \textit{Transportation Science}, 60, 284-315. \url{https://doi.org/10.1287/trsc.2025.0116}

  \item Dongyang Xia, Lixing Yang, Yahan Lu, and Shadi Sharif Azadeh (2025). Robust charging station location and routing-scheduling for electric modular autonomous units. \url{https://arxiv.org/abs/2504.04408}
\end{enumerate}

\newpage

\section{About This Brief}

This brief synthesizes four studies on modularized transit systems, including one charging study using technical data associated with NExT MAVs. It translates the main research findings into deployment implications for MAV manufacturers, public transport operators, cities and regulators, charging partners, autonomous transit software developers, and industrial partners considering pilots. The brief also identifies where data-driven optimization, scenario design, and visualization platforms can support pilot planning and procurement discussions. The quantitative results come from computational experiments and real-world transit-data case studies, and should be interpreted as scenario-specific findings on direction and order of magnitude. Real-world deployment outcomes will depend on local demand, road conditions, station space, charging infrastructure, regulation, vehicle readiness, data availability, and operational integration.
\end{document}